\documentclass[11pt]{article}
\usepackage[margin=1.3in]{geometry}
\usepackage{graphicx} % Required for inserting images
\usepackage{authblk}  % For better author formatting
\usepackage{hyperref} % For clickable email links
\usepackage[dvipsnames]{xcolor}
\usepackage{esvect}
\usepackage[ruled,vlined,linesnumbered]{algorithm2e}
\SetAlgoNoLine
\usepackage{graphicx}      
\usepackage{amsmath,amssymb, amsfonts,amsthm}
\usepackage{subcaption}
\usepackage[ruled,vlined,linesnumbered]{algorithm2e}

\usepackage{algpseudocode}
\usepackage{standalone}
\usepackage{hyperref}   
\usepackage{color}% for hyperlinks
\newtheorem{Theorem}{Theorem}
\newtheorem{Example}{Example}
\newtheorem{Lemma}{Lemma}

\newtheorem{Remark}{Remark}
\newtheorem{proposition}[Theorem]{Proposition} 
\usepackage{tikz}
\usetikzlibrary{arrows.meta,calc}
\usepackage{float}
\usepackage{todonotes}
\usepackage[dvipsnames]{xcolor}
\usepackage{enumitem}
\usepackage{url}
\usepackage{array}
\usepackage{caption}
\usepackage{nicematrix}
\usepackage{arydshln}  
\usepackage{array} % Required for tabular
\usepackage{tikz}
\usepackage{pgfplots}
\pgfplotsset{compat=1.8}
\usepackage{pgfplotstable}
\usetikzlibrary{3d, calc}
\usepackage{placeins}   % Ensures figures do not float past barriers
\usepackage{cases}
\pgfplotsset{compat=1.17}

\usepackage[
    backend = biber,
    doi = false,
    url = false,
    isbn = false,
    sorting = none,
    giveninits = true,  % modern replacement for firstinits
    maxnames = 99,
    hyperref = true,
    block = none
]{biblatex}

\renewbibmacro{in:}{}

\bibliography{reference}

\newtheorem{Definition}{Definition}
\newtheorem{Problem}{Problem}

\newcommand{\N}{\mathbb{N}}

\newcommand{\rank}
{\operatorname{rank}}

\DeclareMathOperator{\diag}{diag}

\newcommand{\R}{\mathbb{R}}

\newcommand{\C}{\mathbb{C}}

\newcommand{\xc}[1]{\vspace{.1cm}
\noindent {\em #1} }

\makeatletter
\def\blfootnote{\xdef\@thefnmark{}\@footnotetext}
\makeatother

\title{Inverse Problem for Super-linearization with Single Observable and Diagonalizable Principal Submatrix}

\begin{document}
\author{ Haoyu Yin$^*$ \quad and \quad Xudong Chen 
}
\date{}
\maketitle
\blfootnote{H. Yin and X. Chen are with Department of Electrical \& Systems Engineering, Washington University in St. Louis. Emails: \texttt{\{h.yin, cxudong\}@wustl.edu}. }

\blfootnote{$^*$Corresponding author: Haoyu Yin.}
%\blfootnote{This work was supported by NSF ECCS-2426017 and by DARPA HR0011-25-3-0211.}

\begin{abstract}                          
Super-linearization is, roughly speaking, embedding an $n$-dimensional nonlinear
dynamical system into an $(n + m)$-dimensional
linear system. 
We introduce and address, in this paper, the inverse problem: 
given a linear system $\dot x= M x$, what type of dynamical systems can be super-linearized to it? 
We focus on the special case where the dimension reduction is one (i.e., $m = 1$) and the leading $n$-by-$n$ principal submatrix of $M$ is diagonalizable, and provide a solution to the problem.

\end{abstract}

% \end{frontmatter}

\section{Introduction}\label{sec:intro}
Super-linearization~\cite{belabbas2022visible,ko2024minimum,ko2025super} is about embedding of a nonlinear dynamical system into a finite, yet higher dimensional linear system. 
Strong super-linearization~\cite{ko2024minimum} further requires that the added coordinates
remain equal to the corresponding observable
functions along the trajectories. %We will present the precise definitions shortly. %by augmenting the state $x(t)$ with a vector of observable functions. 
This notion aligns with classical Koopman embedding and Carleman linearization~\cite{koopman1931hamiltonian,carleman1932application,brunton2016koopman}, the embedding spaces of the latter are however infinite-dimensional.   
Over the past decade, interest in finite-dimensional lifting has surged, driven by its significant computational advantages in prediction and control across diverse applications, such as power system stability assessment~\cite{susuki2016applied}, trajectory tracking for soft robots~\cite{bruder2019modeling}, and modal analysis of fluid flows~\cite{mezic2013analysis}.

Driven by these applications, there has been an active development in understanding when a given nonlinear system admits a super-linearization, how to efficiently check it, what is the canonical form for the associated super-linearized system, etc. The literature review will be provided shortly in Subsection~\ref{subsec:literature_review}.   
In this paper, we investigate super-linearization from a different perspective: instead of looking for super-linearizations of a given nonlinear dynamical system, we start with a linear system of dimension $(n + m)$ and ask: what type of nonlinear systems of dimension~$n$ can be super-linearized to it? 
This problem is new. We focus on the case where $m = 1$ and provide a complete solution to the problem. 
%We now formulate the problem precisely and state the main result of the paper. 

%Constructing such representations is inherently challenging because Lie differentiation of candidate observables along nonlinear vector fields typically introduces higher-order terms, triggering an infinite functional cascade that fails to close within a finite-dimensional linear space~\cite{brunton2016koopman,iacob2023finite}. 

\subsection{Problem Formulation}
In this subsection, we provide a precise problem formulation.  
We start by recalling the definitions of super-linearization.  Consider an autonomous nonlinear system
\begin{equation}
\label{eq:f_x}
\dot{x}=f(x), \qquad x\in\R^n,
\end{equation}
where $f:\R^n\to\R^n$ is a smooth vector field. For each
$x_0\in\R^n$, let $e^{tf}x_0$ be the solution of
\eqref{eq:f_x} at time $t\geq0$ with initial condition
$x(0)=x_0$.

Given a map $p:\R^n\to\R^m$, we augment the state $x(t)\in \R^n$ to the vector $z(t):= \hat p(x(t))$, where
$$
\hat p(x(t))
:= 
\begin{bmatrix}
x(t)\\
p(x(t))
\end{bmatrix}
\in\R^{n+m}.
$$
Following the notation of~\cite{ko2025super}, we define two matrices $\Pi$ and $\bar \Pi$ as follows:
$$
\Pi:=
\begin{bmatrix}
I_n&0
\end{bmatrix}\in \R^{n\times (n + m)},
\qquad
\bar{\Pi}:=
\begin{bmatrix}
0&I_m
\end{bmatrix} \in \R^{m\times (n + m)}.
$$
We further have the following definition:

\begin{Definition}
\label{def:superlinearization}
A pair $(M,p)$, with
$M\in\R^{(n+m)\times(n+m)}$ and $p:\R^n\to \R^m$, is a
{\bf super-linearization} of~\eqref{eq:f_x} (or simply the vector field~$f$) if
\begin{equation}
\label{eq:super}
\Pi z(t)=x(t), \quad \mbox{with } z(t) := e^{tM}\hat p(x_0)  \mbox{ and } x(t):= e^{tf}x_0, 
\end{equation}
for every $x_0\in\R^n$ and every $t\geq0$ for which the
nonlinear solution exists \cite{belabbas2022visible}. 
The super-linearization is {\bf{strong}} if, in addition, for any $t\geq 0$,
\begin{equation}
\label{eq:proj_strong}
    \bar{\Pi}z(t)
    =
    p(x(t)).
\end{equation}
The entries of $p$ are called the {\bf observables}. 

% The last $m$ coordinates must therefore remain equal to
% the observables evaluated along the trajectory.
\end{Definition}

In words, condition~\eqref{eq:super} requires that the first $n$ coordinates of the solution $z(t)$ of the linear system $\dot z(t) = Mz(t)$ coincide with the solution $x(t)$ of the nonlinear system $\dot x(t) = f(x(t))$.
%$e^{tM}z_0$ reproduce the state trajectory. 
Condition~\eqref{eq:proj_strong} requires that the last $m$ coordinates of the solution $z(t)$ be identical with $p(x(t))$.  Under strong super-linearization, the original state coordinates $x_1, \ldots, x_n$ together with the observables $p_1, \ldots, p_m$ span a finite-dimensional space that is invariant under the dynamical system~\eqref{eq:f_x}.

Let $(M,p)$ be a strong super-linearization of~\eqref{eq:f_x}. 
We partition the matrix $M \in \mathbb{R}^{(n+m) \times(n+m)}$ into $2$-by-$2$ matrix blocks as follows: 
\begin{equation}
\label{eq:def_M}
    M = \begin{bmatrix} A & b \\ c^\top & d \end{bmatrix},
\end{equation}
where  $A \in \mathbb{R}^{n \times n}$, $b \in \mathbb{R}^{n\times m}$, $c \in \mathbb{R}^{n \times m}$, and $d \in \mathbb{R}^{m \times m}$. 
%The linear system $\dot{z}(t)=M z(t)$ generates a strong %superlinearization for $\dot{x}(t) = f(x(t))$ if and only if 
It directly follows from~\eqref{eq:super} and~\eqref{eq:proj_strong} that the vector field $f$ and the map $p$ satisfy 
\begin{align}
    f(x)&=Ax+bp(x), \label{eq:f_invariance} \\
    \frac{\partial p(x)}{\partial x} (Ax+bp(x))&=c^\top x+dp(x).
        \label{eq:invariance}
\end{align}
Conversely, if~\eqref{eq:f_invariance} and~\eqref{eq:invariance} are satisfied, then $(M,p)$ is a strong super-linearization of~\eqref{eq:f_x}. For a later purpose, we introduce the following definition: 

\begin{Definition}
    Given the matrix $M$ in~\eqref{eq:def_M}, a function $p: \R^n\to \R^m$ is {\bf compatible with} $M$ if it satisfies~\eqref{eq:invariance}. 
\end{Definition}

It is clear that every compatible $p:\R^n\to \R^m$ gives rise to a vector field~$f:\R^n\to \R^n$ via~\eqref{eq:f_invariance}. 
The problem we will aim to address is the following:

\begin{Problem}\label{prob:main}
    Given an $(n + m)$-by-$(n + m)$ matrix $M$, find 
    all compatible, polynomial maps $p:\R^n\to \R^m$. 
    %all nonlinear {\color{blue} polynomial} vector fields $f: \R^n\to \R^n$, together with the {\color{blue} polynomial} maps $p:\R^n\to \R^m$ that is compatible with $M$, such that $(M,p)$ is a strong superlinearization of $f$.
\end{Problem}

In this work, we take a first step toward Problem~\ref{prob:main} by considering the minimal state-retaining augmentation with a scalar polynomial observable ($m = 1$), lifting the system state from $\R^n$ to $\mathbb{R}^{n+1}$. %A single observable is the smallest
% state-retaining augmentation capable of
% representing nonlinear dynamics, as setting $m = 0$ restricts $f(x)$ to be purely linear.
% We focus on state matrices $A$ that are diagonalizable. This class encompasses standard benchmark constructions such as the triangular monomial system in \cite[Eqs.~(21)--(24)]{brunton2016koopman}, as well as real symmetric state blocks, including negative graph Laplacians governing continuous-time consensus dynamics on undirected networks~\cite{olfati2004consensus}.
Under the assumption that the submatrix~$A$ of $M$ is diagonalizable (over $\C$), we provide in Theorem~\ref{thm:main} a complete characterization of polynomial maps $p:\R^n\to \R$ that are compatible. 

\subsection{Motivating Examples}

%In contrast to traditional forward approaches that construct linear embeddings for a given nonlinear vector field $f(x)$, we investigate strong superlinearization from an inverse perspective in this work. 

%Given a prescribed evolution matrix $M$, a natural question is how to characterize all observable functions $p(x)$ such that the system satisfies~\eqref{eq:inv} and can be superlinearized into the prescribed linear system. 
The following examples illustrate that polynomial observables
need not exist for a prescribed lifted matrix. Throughout the examples,
we fix the ordering of the coordinates. The first $n=3-m$ coordinates
are retained as the states and the remaining $m$ coordinates
are the observables.

\begin{Example}[Nonexistence for both partitions]
\label{ex:nonexistence}
\normalfont
Consider
$$
M_1=\begin{bmatrix}
-1&0&0\\0&-1&1\\1&0&-1
\end{bmatrix}.
$$

\xc{Case 1: $m=1$.}
In this case, we have that $A=-I_2$, $b=(0,1)^\top$, $c=(1,0)^\top$,
and $d=-1$. Write $p(x)=e+\ell^\top x+h(x)$, where $h$
contains only terms of degree at least two. We will see later in Theorem~\ref{thm:main} that compatibility of~$p$ requires that
$$
\alpha=-1-\ell_2,\qquad
c=A^\top\ell-\alpha\ell
=\begin{bmatrix}\ell_1\ell_2\\\ell_2^2\end{bmatrix}.
$$
Thus, $\ell_1\ell_2=1$ and $\ell_2^2=0$, which is a contradiction.

\xc{Case 2: $m=2$.}
The blocks are $A=-1$, $b=(0,0)$, $c=(0,1)$, and
$d=\bigl[\begin{smallmatrix}-1&1\\0&-1\end{smallmatrix}\bigr]$.
Hence $f(x)=Ax+bp(x)=-x$. For $p=(p_1,p_2)^\top$, the second
component of the vector-valued form of~\eqref{eq:invariance} gives
\begin{equation}
\label{eq:examples-linear-obstruction}
x\frac{\partial p_2}{\partial x}-p_2=-x.
\end{equation}
For every polynomial $p_2$, the coefficient of $x$ on the left is
zero, whereas that on the right is $-1$. Thus no polynomial
observable map exists for this partition either.
\end{Example}

\begin{Example}[Existence for $m=1$ but not for $m=2$]
\label{ex:existence-m1}
\normalfont
Consider
$$
M_2=\begin{bmatrix}
-1&0&0\\0&-2&1\\1&0&-1
\end{bmatrix}.
$$

\xc{Case 1: $m=1$.}
Now we have that $A=\operatorname{diag}(-1,-2)$, with the same $b,c,d$
as in the first example. We establish later in Theorem~\ref{thm:main} that the observable satisfies
$$
p(x_1,x_2)=x_1+x_2+a x_1^2,\qquad a\in\R.
$$
Indeed, the affine conditions give $\ell=(1,1)^\top$, $e=0$,
and $\alpha=-2$, and the condition on the nonlinear components gives $h=a x_1^2$.
Direct substitution into
\eqref{eq:f_invariance}--\eqref{eq:invariance} verifies
\begin{align*}
f(x)&=Ax+bp(x)
=\begin{bmatrix}-x_1\\x_1-x_2+a x_1^2\end{bmatrix},\\
\frac{\partial p}{\partial x}(Ax+bp(x))
&=-x_2-a x_1^2=c^\top x+dp(x).
\end{align*}
The state dynamics are nonlinear for $a\ne0$.

\xc{Case 2: $m=2$.}
The blocks $A,b,c$ are unchanged from the case where $m=2$ in Example~\ref{ex:nonexistence}, and
$d=\bigl[\begin{smallmatrix}-2&1\\0&-1\end{smallmatrix}\bigr]$.
The second observable equation is again
\eqref{eq:examples-linear-obstruction}, so no polynomial
observable map exists for this partition.
\end{Example}

\subsection{Literature Review}
\label{subsec:literature_review}

Existing studies on exact finite-dimensional
linear representations mainly adopt a forward
perspective: given a nonlinear vector field
$f(x)$, they seek a lifted matrix $M$
and observables $p$ such that the evolution
admits an exact finite-dimensional linear
representation. Research in this direction
has focused on establishing existence
conditions and constructing compatible
observables through data-driven and
analytical methods. %In contrast, the inverse
% problem of characterizing all observables
% compatible with a prescribed lifted matrix
% $M$ remains largely open.

% Below, we review the literature most closely related to our work.

Data-driven methods typically search for
Koopman-invariant subspaces within a prescribed
set of candidate functions. For instance, the symmetric subspace decomposition
method in~\cite{haseli2021learning} identifies
the maximal Koopman-invariant subspace contained
in the span of a given dictionary under independent and identically distributed sampling assumptions. The resulting characterization is therefore
restricted to the supplied candidates.  
For more general polynomial systems, the graph-theoretic approach
in~\cite{belabbas2023sufficient} establishes
a sufficient condition for polynomial
super-linearizability based on constant products
of edge weights along directed cycles
of a weighted dependency graph.
A procedure for
constructing the corresponding observables has been provided by the authors.

A closely related line of research investigates the
structure of super-linearizations and the
number of observables required. The work in~\cite{belabbas2022visible}
distinguishes visible and hidden observables
according to whether they appear explicitly
in the original state dynamics and derives
a tight lower bound on the number of
visible observables.
For polynomial systems with an equilibrium
at the origin that admit balanced
super-linearizations with a single visible
observable, the work
in~\cite{belabbas2022canonical} establishes
a canonical form through a linear change
of state variables. The resulting form separates an autonomous
linear subsystem from the nonlinear dynamics,
whose nonlinear terms depend only on the
state of the linear subsystem. The work in~\cite{ko2024minimum} further
identifies system invariants associated with
super-linearizations and determines the
minimum dimension of a linear representation
when one exists.
More recently,~\cite{ko2025super}
establishes necessary and sufficient conditions
for polynomial systems with a single visible
monomial observable to admit strong
super-linearizations using monomial observables.
The work also shows that the additional
observables essential to the representation
have the same polynomial degree as the
visible monomial.

The existence of finite-dimensional linear
representations has also been investigated
from geometric and algebraic perspectives.
The work in~\cite{kvalheim2026linearizability}
establishes necessary and sufficient conditions
for global linearizing embeddings of dynamical
systems on compact state spaces or with
compact attractors.
For polynomial systems, the classical works
in~\cite{coomes1991linearization,van1994locally}
establish relationships between polynomial
flows and locally finite polynomial derivations,
providing an algebraic perspective on the
existence of exact polynomial representations.

\noindent\textbf{Notations.}
We summarize the notation used throughout the paper.
We denote by $\R$, $\C$, and $\N_0$ the sets
of real numbers, complex numbers, and nonnegative
integers, respectively. The matrix $I_n$ denotes
the $n\times n$ identity matrix.
For a polynomial $h(x)$, we denote by
$h^{(m)}(x)$ its homogeneous component of
degree $m$.

For any $\beta\in\N_0^q$ and $y\in\C^q$,
we define
$|\beta|:=\sum_{i=1}^q\beta_i,
y^\beta:=\prod_{i=1}^q y_i^{\beta_i}.
$
% In particular, for $\eta^{(i)}\in\N_0^q$,
% we have that $
% y^{\eta^{(i)}}
% =\prod_{j=1}^q y_j^{\eta_j^{(i)}}$.
For any $\eta,\gamma\in\N_0^q$, we write
$\gamma\leq\eta$ if $\gamma_i\leq\eta_i$
for every $i$. Throughout this paper, kernels and column spaces
involving complex matrices are understood over
$\C$, unless otherwise stated.

\section{Main Result}
\label{sec:charc}
We need a few preliminaries to state the main result of the paper. 
Given $(A, b)\in \R^{n\times n}\times \R^n$, let $\mathcal{C}(A, b)$ be the controllability matrix: 
$$ \mathcal{C}(A,b):=
\begin{bmatrix}
b & Ab & \cdots & A^{n-1}b
\end{bmatrix},$$ 
and 
$$r:=\operatorname{rank} \mathcal{C}(A,b) \quad \mbox{and} \quad q:=n-r.$$
Let $\mu=(\mu_1,\mu_2,\ldots,\mu_s)$, where $\mu_1,\mu_2,\ldots,\mu_s$ are distinct eigenvalues of $A$, and 
$$
    E_a:=\ker (A-\mu_a I_n) \quad \mbox{and} \quad q_a:=\dim E_a.
$$
We assume in the sequel that $A$ is diagonalizable over $\C$.
It follows that $\C^n$ is a direct sum of $E_a$, for $a = 1,\ldots, s$. 
Thus, there exists a unique decomposition
$$b=\sum_{a=1}^sb^{(a)},\quad \mbox{with }  b^{(a)}\in E_a.$$ 
Using the vectors $b^{(a)}$, we can express the column space of $\mathcal{C}(A, b)$ as
\begin{equation}\label{eq:imC}
\operatorname{col}\mathcal{C}(A, b)=\operatorname{span}\left\{b^{(a)}: b^{(a)} \neq 0\right\}.
\end{equation}

% {\color{blue} We use $\operatorname{col}_{\mathbb{R}} \mathcal{C}(A, b)$ and $\operatorname{col}_{\mathbb{C}} \mathcal{C}(A, b)$ to denote the real and complex column spans, respectively.} 

% {\color{blue} from $\mu$ to $\lambda$ there is a gap}

For $q>0$, the following result follows from the controllability
decomposition~\cite[Theorem~6.6]{chen2013linear}
and the diagonalizability of $A$ (note that~\cite[Theorem~6.6]{chen2013linear} is stated for a real matrix $Q$, but the arguments can be extended straightforwardly to the complex case).
%A proof is provided in the appendix.

% The following lemma allows us to exploit the structure of the observable $p(x)$ induced by the diagonalizable matrix $A$.
% \begin{Lemma}
% \label{lem:imC}
%     Let $A$ be diagonalizable over $\R$. Then $$\operatorname{col}\mathcal{C}(A, b)=\operatorname{span}\left\{b^{(a)}: b^{(a)} \neq 0\right\}.$$
% \end{Lemma}

\begin{Lemma}
\label{lem:ker_Q}
    %Let $A$ be diagonalizable over $\R$. Then 
    For $q > 0$, there exists a full row-rank matrix $Q\in \C^{q\times n}$ and a vector $\lambda\in \C^q$ 
    such that
    \begin{equation}
    \label{eq:def_Q}
        \ker Q=\operatorname{col} \mathcal{C}(A,b), \quad Qb=0,\quad \mbox{and} \quad QA=\diag(\lambda) Q.
    \end{equation}
\end{Lemma}

%Lemma~\ref{lem:ker_Q} constructs projected coordinates that eliminate the subspace while preserving a diagonal representation of $A$. The construction is given in \cite{}. We note that such $Q$ is not unique given different choices of basis for the repeated eigenvalue coordinates $b^{(a)}$. 

We pick any such matrix $Q\in \C^{q\times n}$, and define 
\begin{equation}
\label{eq:new_coordinate}
    y:=Qx \in \C^q.
\end{equation}
%For a vector $\eta\in \N_0^q$, we let $|\eta|:= \sum_{i = 1}^q \eta_i$. 
Given $\alpha\in \R$ and an integer $k\in \N_0$, we let
\begin{equation}
\label{eq:S_k}
\mathcal{S}^k_\alpha :=\left\{\eta \in \mathbb{N}_0^q: |\eta| \geq k \mbox{ and } \eta^\top \lambda =\alpha  \right\}. 
\end{equation}
In case $k = 0$, we simply write $\mathcal{S}_\alpha$. 
We will mainly focus on the sets $\mathcal{S}_0$ and $\mathcal{S}^2_\alpha$.
% An arbitrary convex cone $\mathcal{S}$ is said to be {\it a rational polyhedral cone} (see, e.g.,~\cite{}) if it takes the form $\mathcal{S} := \{x : A x \geq b \}$
% there exist a finite number of generators such that every element in the cone can be expressed as a rational, nonnegative combination of the generators.    

For $\kappa,\eta\in\N_0^q$, we write
$\kappa\leq\eta$ if $\kappa_i\leq\eta_i$
for every $i$. An element $\kappa \in \mathcal{S}^k_\alpha$ is said to be {\it component-wise minimal} (or simply {\it minimal}) if there is no distinct $\eta \in \mathcal{S}^k_\alpha$ satisfying $\eta \leq \kappa$.
We need the following result: 

\begin{Lemma}
\label{lem:finite_generators}
There exists a finite number of generators $\eta^{(1)}, \ldots, \eta^{(K)} \in \mathcal{S}_0$,
such that every element in $\mathcal{S}_0$ can be expressed as a linear combination of the generators with nonnegative integer coefficients.
\end{Lemma}

\begin{proof}
By Dickson's lemma,
$\mathcal S_0-\{0\}$ has finitely
many minimal elements,
denoted by
$\eta^{(1)},\ldots,\eta^{(K)}$. 
We claim that every element $\eta\in \mathcal{S}_0$ can be written as $\eta = \sum_{j=1}^K a_j \eta^{(j)}$, where $a_j\in\N_0$.
We prove this claim by induction on $|\eta|$. 
For the base case where $\eta = 0$, the claim holds trivially. 
For the inductive step, we assume that the claim holds for any $\xi$ with $|\xi| \leq k$, and we prove for $\eta$ with $|\eta| = k + 1$. 
% For any nonzero $\eta\in\mathcal S_0$,
% suppose that every $\xi\in\mathcal S_0$
% with $|\xi|<|\eta|$ can be expressed as a linear combination of the generators with non-negative integer coefficients,
% we show that $\eta$ also admits such a
% representation. 
Since $\eta\in\mathcal S_0-\{0\}$,
there exists a minimal element
$\eta^{(i)}\leq\eta$. Let $\xi:=\eta-\eta^{(i)}$, we have that
$\xi\in\N_0^q$ and $\lambda^\top\xi=0$. 
Moreover,
since $\eta^{(i)}\neq0$, we have
$|\xi|<|\eta|$.
By the induction hypothesis, there exist
$a_1,\ldots,a_K\in\N_0$ such that
$\xi=\sum_{j=1}^K a_j\eta^{(j)}.
$
We can write
$\eta=\eta^{(i)}+\xi
=\eta^{(i)}+\sum_{j=1}^K a_j\eta^{(j)}$, 
which is a nonnegative integer combination
of the minimal elements.
This completes the proof.
\end{proof}
% The lemma follows from Dickson's Lemma~\cite[Chapter 2, Section 4]{cox2007ideals}
% and the additive structure of $\mathcal S_0$.
% A proof is provided in
% Appendix~\ref{app:finite_generators}. 

By Dickson's Lemma, we have that $\mathcal{S}^2_\alpha$ has finitely many minimal elements, denoted by $\kappa^{(1)},\kappa^{(2)},\ldots, \kappa^{(N)}$. We are now ready to state the main result of the paper:
%We then introduce 
%$$g_i(y):=y^{\eta^{(i)}}, i=1,\ldots,K.$$

% We further define the sets
% \begin{equation}
% \label{eq:semigroup}
% \begin{aligned}
% \mathcal{S}_0&:=\left\{\eta \in \mathbb{N}_0^q: \lambda^{\top} \eta=0\right\},  \\ 
% \mathcal{S}_\alpha^{\geq 2}&:=\left\{\beta \in \mathbb{N}_0^q:|\beta| \geq 2, \lambda^{\top} \beta=\alpha\right\},
% \end{aligned}
% \end{equation}
%where $|\eta|:=\sum_{i=1}^q \eta_i$. 

\begin{Theorem}
\label{thm:main}
    Let $M = [A, b; c^\top, d ]\in \R^{(n+1)\times (n + 1)}$, with $A\in \R^{n\times n}$ diagonalizable over $\C$. Let $p: \R^n\to \R$ be a real polynomial in $x\in \R^n$ and we write
    \begin{equation}
    \label{eq:px}
    p(x) = e + \ell^\top x + h(x), 
    \end{equation}
    where $e\in \R$, $\ell\in \R^n$, and $h$ is a polynomial with no nonzero affine terms. Then, $p$ is compatible with~$M$ if and only if the following two items hold:
    \begin{enumerate}
    \item Let $\alpha: = d - \ell^\top b$. Then, \begin{subequations}\label{eq:constraints}
        \begin{align}
            \alpha e &= 0,
        \label{eq:constraint3}\\
        c-A^\top \ell+\alpha\ell &= 0,
        \label{eq:constraint4}
        \end{align}
        \end{subequations}
    \item 
    Let $q:= n - \operatorname{rank}\mathcal{C}(A, b)$. Then, there are two cases: 
    \begin{enumerate}
        \item If $q = 0$, then $h = 0$.
        \item If $q > 0$, then we let $Q\in \C^{q\times n}$ and $\lambda\in \C^{q}$ be such that~\eqref{eq:def_Q} is satisfied. Let $\eta^{(1)},\eta^{(2)},\ldots,\eta^{(K)}$ be the generators of $\mathcal{S}_0$, and $\kappa^{(1)},\kappa^{(2)},\ldots,\kappa^{(N)}$ be the minimal elements of $\mathcal{S}_\alpha^2$. Then, $h$ takes the following form:
        \begin{equation}
        \label{eq:h_form}h(x)=\operatorname{Re}\left[\sum_{j=1}^N y^{\kappa^{(j)}} F_j\left(y^{\eta^{(1)}}, \ldots, y^{\eta^{(K)}}\right)\right], \quad \mbox{with }y=Qx,
        \end{equation}
        where $F_j$ are arbitrary complex polynomials in $K$ complex variables.
    \end{enumerate}
    \end{enumerate}
\end{Theorem}

%  Let $\alpha:=d-\ell^\top b$, Then, a polynomial $p:\R^n\rightarrow \R$ satisfies~\eqref{eq:invariance} if and only if 
%     $p$ takes the form: 
% \begin{equation}
% \label{eq:p_form}
%     p(x)=e+\ell^{\top} x+\operatorname{Re}\left[\sum_{j=1}^N y^{\kappa^{(j)}} F_j\left(y^{\eta^{(1)}}, \ldots, y^{\eta^{(K)}}\right)\right], \mbox{ with }y=Qx,
% \end{equation}
% where $\ell\in \R^n$ and $e\in \R$ satisfy 
% \begin{subequations}\label{eq:constraints}
% \begin{align}
% (d-\ell^\top b) e &= 0,
% \label{eq:constraint3}\\
% c-A^\top \ell+\alpha\ell &= 0,
% \label{eq:constraint4}
% \end{align}
% \end{subequations}
% and $F_j$ are arbitrary complex polynomials.

We have the following two remarks:

\begin{Remark}
\label{rem:alpha}
\normalfont
% By Theorem~\ref{thm:main}, the polynomial observables $p$ can be constructed via the following procedure.
% First, 
Note that~\eqref{eq:constraint4} is equivalent to the following equation: 
$$
M^\top
\begin{bmatrix}
-\ell\\1
\end{bmatrix}
=
\alpha
\begin{bmatrix}
-\ell\\1
\end{bmatrix}.
$$
%In words, for $\alpha$ and $\ell$ to satisfy~\eqref{eq:constraint4},  $\alpha$ must be a real eigenvalue
%of $M^\top$, and $(-\ell, 1)$ must be the corresponding eigenvector. 
Thus, to obtain a compatible polynomial~$p$ for a given matrix $M$, we can first find a real eigenvalue
$\alpha$ of $M^\top$, if it exists, whose eigenspace
contains a real eigenvector $v$ with last entry being nonzero.  
Normalize the eigenvector~$v$
so that $v = (-\ell; 1)$. In this way, we obtain
the vector $\ell$ that satisfies the condition~\eqref{eq:constraint4}.
Next, we choose $e$ so that~\eqref{eq:constraint3} is satisfied. Specifically, $e$ can be arbitrary if $\alpha=0$ and $e = 0$ otherwise. Then, we compute the rank $r$ of
$\mathcal C(A,b)$ and let $q:=n-r$.
If $q=0$, we set $h(x)=0$ and $p(x)=e+\ell^\top x$.
Otherwise, we obtain $Q$ and $\lambda$
as in Lemma~\ref{lem:ker_Q}.
With $\alpha$ and $\lambda$ given, we compute the $K$ generators of
$\mathcal S_0$ and the $N$ minimal elements of $\mathcal S_\alpha^2$; the former can be obtained by using the algorithm in~\cite{bruns2010normaliz} and the latter by the algorithm in~\cite{pottier1991minimal}. 
Finally, choosing a set of complex
polynomials $F_1,\ldots,F_N$ in $K$ variables, 
we obtain $h$ as in~\eqref{eq:h_form} and
hence, $p(x) = e + \ell^\top x + h(x)$. \hfill{\qed}
\end{Remark}

\begin{Remark}
\label{rem:dependecy}
\normalfont
The family of real polynomial observables
characterized by Theorem~\ref{thm:main}
is independent of the choice of $Q$
satisfying~\eqref{eq:def_Q} and the
generators of $\mathcal S_0$.
Different choices of $Q$ may yield different
coordinates $y=Qx$, while different
generating sets of $\mathcal S_0$ may change
the monomials $y^{\eta^{(i)}}$ and the
corresponding polynomials $F_j$
in~\eqref{eq:h_form}.
Nevertheless, the resulting representations
characterize the same family of admissible
real polynomials $p(x)$. \hfill{\qed}
\end{Remark}

%\subsection{Illustration of Theorem~\ref{thm:main}}

Toward the end of the section, we provide an example to illustrate Theorem~\ref{thm:main} and Remark~\ref{rem:alpha}.

\begin{Example}
\label{ex:main}
\normalfont
Consider the matrix $M = [A, b; c^\top, d]$, with
$$A =\begin{bmatrix}
    1 & -1& 0 & 0 & 0\\
    1 &  1& 0 & 0 & -1\\
    0 &  0& -1 & -2 & 2\\
    0 &  0&  2& -1 & -2\\
    0 &  0& 0 & 0 & 1
\end{bmatrix}, \quad b=\begin{bmatrix}
    1\\ 0\\1\\0\\1
\end{bmatrix},\quad c=\begin{bmatrix}
    0\\0\\0\\0\\-1
\end{bmatrix}, \quad d=3.$$
The set of eigenvalues of $A$ is given by $\{1,1-i,1+i,-1+2i,-1-2i\}$, which are pairwise distinct. Thus, $A$ is diagonalizable over $\C$, which satisfies the hypothesis of Theorem~\ref{thm:main}. Moreover, we have that $r=\rank \mathcal{C}(A,b)=1$ and $q=n-r=4$. 

To obtain a compatible polynomial~$p$, we first compute the eigenvalues of $M^\top$. In this case, $M$ has a unique real eigenvalue $\alpha=2$ whose eigenspace is spanned by $(0,0,0,0,-1,1)$, so we set $\ell:=(0,0,0,0,1)^\top$.  Since $\alpha\neq 0$, we must have that $e=0$. 
We then define
$$Q:=\begin{bmatrix}
    1& i& 0& 0 &-1\\
    1& -i& 0& 0 &-1\\
    0&  0& 1& i &-1\\
    0&  0& 1& -i &-1
\end{bmatrix} 
\quad \mbox{and} \quad 
\lambda:=(1+i,1-i,-1+2i,-1-2i),$$ 
which satisfy~\eqref{eq:def_Q}. 
The entries of $y=Qx$ are explicitly given by
$$y_1=x_1-x_5+ix_2,\quad y_2=x_1-x_5-ix_2,\quad y_3=x_3-x_5+ix_4,\quad y_4=x_3-x_5-ix_4.$$
Using the algorithm in~\cite{bruns2010normaliz}, we obtain the following elements that generate $\mathcal{S}_0$:
$$\eta^{(1)}=(1,1,1,1), \quad \eta^{(2)}=(4,0,1,3),\quad  \eta^{(3)}=(0,4,3,1).$$
Further, using the algorithm in~\cite{pottier1991minimal}, we obtain the minimal elements of $\mathcal{S}_\alpha^2$, for $\alpha = 2$, as follows: 
$$\kappa^{(1)}=(1,1,0,0), \quad \kappa^{(2)}=(4,0,0,2),\quad  \kappa^{(3)}=(0,4,2,0).$$
%Theorem~\ref{thm:main} then characterizes the complete family of observables $p$ as
Then, any compatible polynomial takes the following form:
$$p(x)=x_5+\operatorname{Re}\left[\sum_{j=1}^3 y^{\kappa^{(j)}} F_j\left(y_1y_2y_3y_4,y_1^4y_3y_4^3,y_2^4y_3^3y_4\right)\right],$$
where $F_1, F_2, F_3$ are arbitrary complex polynomials.

\end{Example}

\section{Proof of Theorem~\ref{thm:main}}
\label{sec:compatibility}
% Let $M$ be given as in~\eqref{eq:def_M}.
% For a polynomial function $p(x)$,  we decompose $p$ as follows: \begin{equation}
%     \label{eq:px}
%         p(x)=e+\ell^\top x+h(x),
% \end{equation} 
% where $e\in\mathbb{R}$, $\ell\in\mathbb{R}^{n}$, and $h(x)=\sum_{m=2}^k h^{(m)}(x)$, with $h^{(m)}(x)$ the homogeneous polynomial of degree~$m$. We further define
% $$
% \alpha: = d - \ell^\top b,
% $$
% which plays an important role in the analysis, and helps ease the notation.  
% We substitute $p$ in~\eqref{eq:invariance} with the decomposition~\eqref{eq:px}. The following lemma then states a necessary and sufficient condition on $(e,\ell,h)$ for $p$ to satisfy~\eqref{eq:invariance}.
Let $M$, $p$ and $\alpha$ be given as in
Theorem~\ref{thm:main}.
%We first substitute~\eqref{eq:px}
%into~\eqref{eq:invariance}.
We first establish the necessary and sufficient
conditions on $(e,\ell,h)$ for $p$ to
satisfy~\eqref{eq:invariance}.

\begin{Lemma}
\label{lem:matching}
    %Let $p$ be defined as in~\eqref{eq:px}. 
    The polynomial $p$ satisfies~\eqref{eq:invariance} if and only if $(e,\ell)$ satisfies~\eqref{eq:constraints} 
    % \begin{align}
    % \label{eq:constraint3}
    % \alpha e&=0,\\ 
    % \label{eq:constraint4}
    % C&=A^\top \ell-\alpha \ell
    % \end{align}
    and, moreover, $(h,\alpha)$ satisfies %for every $m=2,3,\ldots,k$,
    % \begin{align}
    % \label{eq:constraint1}
    %     \nabla h(x)^{\top} b&=0, \\ 
    %     \label{eq:constraint2}
    %     \nabla h(x)^{\top} A x&=\alpha h(x).
    % \end{align}
\begin{subequations}
\label{eq:cond}
\begin{align}
\frac{\partial h(x)}{\partial x} b&=0,
\label{eq:constraint1}\\
\frac{\partial h(x)}{\partial x} A x&=\alpha h(x),
\label{eq:constraint2}
\end{align}
\end{subequations}
for all $x\in \R^n$. 
% \begin{aligned}
    % \label{eq:constraint1}
    %     \nabla h(x)^{\top} b&=0, \\ 
    %     \label{eq:constraint2}
    %     \nabla h(x)^{\top} A x&=\alpha h(x).
    % \end{aligned}
% \end{equation}
\end{Lemma}

Lemma~\ref{lem:matching} follows by substituting $p$ in~\eqref{eq:invariance} with the decomposition~\eqref{eq:px} and matching terms
of the same total degree of the associated homogeneous polynomials. %$h^{(m)}(x)$. Then~\eqref{eq:cond} follows from the preservation of the degree. 
We provide a proof in Appendix~\ref{app:matching_proofs}.

% By Lemma~\ref{lem:matching}, the nonlinear
% component $h$ of an admissible polynomial
% observable must satisfy

% Since the above conditions preserve homogeneous
% degree, and are equivalent to~\eqref{eq:constraint1} and~\eqref{eq:constraint2}.

% In this section, we will establish a
% first derive necessary and sufficient
% conditions on $e$, $\ell$, and $h$ for $p$ to
% satisfy~\eqref{eq:invariance}, so that $p$ is compatible with $M$.

Now, assuming that~\eqref{eq:constraints} holds, we show below that for $p$ to satisfy~\eqref{eq:invariance}, it is necessary and sufficient for its nonlinear part to satisfy $h=0$ when $q=0$, and $h$ to take the form of~\eqref{eq:h_form} when $q>0$. 
% \begin{equation}
% \label{eq:h_form}
% h(x)=\operatorname{Re}\left[\sum_{j=1}^N y^{\kappa^{(j)}} F_j\left(y^{\eta^{(1)}}, \ldots, y^{\eta^{(K)}}\right)\right], \mbox{ with } y=Qx.
% \end{equation}
We take a two-step approach
to prove this, as we outline below. 
\begin{enumerate}
\item 
First, Proposition~\ref{prop:y_alpha}
shows that $h$ satisfies~\eqref{eq:cond}
if and only if $h=0$ when $q=0$.
For $q>0$, it characterizes all nonlinear polynomials
$h$ satisfying~\eqref{eq:cond} using
the admissible exponents in $\mathcal S_\alpha^2$. Specifically, we show that every such polynomial $h$ takes the form $h(x)= \operatorname{Re}\left[\sum_{\beta \in \mathcal{T}}c_\beta (Qx)^{\beta}\right]$, where $\mathcal{T}$ is a finite subset of $\mathcal{S}^2_\alpha$, and vice versa. 
\item Then, in Proposition~\ref{prop:finite}, we show that any complex polynomial $\hat h(y)$ is in the form of $\sum_{\beta\in \mathcal{T}}c_\beta y^\beta$ if and only if it can be expressed as 
$$\hat h(y) = \sum_{j=1}^N y^{\kappa^{(j)}} F_j\left(y^{\eta^{(1)}}, \ldots, y^{\eta^{(K)}}\right).$$
%using the generators of $\mathcal S_0$ and
%the componentwise minimal elements of
%$\mathcal S_\alpha^2$.
\end{enumerate}
%We take the real parts of each complex polynomial $\hat h(y)$

For the remainder of the section, we fix $\alpha\in \R$. %and investigate when a nonlinear polynomial function $h$ can satisfy~\eqref{eq:cond}. 

\subsection{Step~1: Characterization of~$h$}
We first state the following proposition: 
% The following proposition provides an explicit
% representation of the nonlinear components
% $h$.
\begin{proposition}
\label{prop:y_alpha}
The following two items hold:
\begin{enumerate}
\item If $q = 0$, then a real polynomial $h$
with no nonzero affine terms
satisfies~\eqref{eq:cond} if and only if 
$h=0$. 
\item For $q>0$, we let $Q$, $\lambda$
be given as in~\eqref{eq:def_Q} and $\mathcal{S}^2_\alpha$ be given as in~\eqref{eq:S_k}.  
Then, a real polynomial $h$
with no nonzero affine terms
satisfies~\eqref{eq:cond} if and only if there exist a finite subset $\mathcal{T} \subseteq \mathcal{S}^2_\alpha$ and complex numbers $c_\beta$, for $\beta \in \mathcal{T}$, such that
\begin{equation}
\label{eq:h_spectral}
h(x)= \operatorname{Re}\left[\sum_{\beta \in \mathcal{T}}c_\beta (Qx)^{\beta}\right], 
% \left[
% \sum_{\beta\in\mathcal S_\alpha^2}
% c_\beta(Qx)^\beta
% \right],
\end{equation}
for all $x\in \R^n$. 
%where $c_\beta\in\C$ are arbitrary and the sum is finite.
\end{enumerate}
\end{proposition}

% {\color{blue} double check Re in the proof}

To prove Proposition~\ref{prop:y_alpha}, we need a few supporting lemmas.
The following lemma identifies the subspace on which
every polynomial satisfying~\eqref{eq:cond}
is constant.

\begin{Lemma}
\label{lem:iff}
 Let $h$ be a polynomial that satisfies~\eqref{eq:cond}.
    % \begin{equation}
    % \label{eq:cond}
    %     \nabla h(x)^\top b=0,\quad \nabla h(x)^\top Ax=\alpha h(x), \mbox{ for any } x\in \R^n.
    % \end{equation}
    Then,  
    \begin{equation}
    \label{eq:cond_controllability}
    \frac{\partial h(x)}{\partial x}  v=0 \mbox{ for every }v\in \operatorname{col}_\R \mathcal{C}(A,b), \ x\in\R^n. 
    \end{equation}
\end{Lemma}

\begin{proof}
We first prove by induction that for any $j\geq 0$, the following holds
\begin{equation}
    \label{eq:induction_assump}
\frac{\partial h(x)}{\partial x} A^j b=0 \quad \mbox{ for every }x \in \mathbb{R}^n.
\end{equation}
For the case where $j=0$, $\frac{\partial h(x)}{\partial x} b=0$ follows from~\eqref{eq:cond}. We now assume that for some $j\geq 0$,~\eqref{eq:induction_assump} holds
and show that $\frac{\partial h(x)}{\partial x} A^{j+1} b=0$.
First, by the induction hypothesis,
we have that
\begin{equation}
\label{eq:hessian_Ax}
\frac{\partial}{\partial x}
\left(
\frac{\partial h(x)}{\partial x}A^jb
\right)
=
(A^jb)^\top
\frac{\partial^2h(x)}{\partial x^2}
=0.
\end{equation}
It follows that
$$
(Ax)^\top
\frac{\partial^2h(x)}{\partial x^2}
A^jb=0,
$$
where we have used the symmetry
of the Hessian matrix.

Next, using~\eqref{eq:constraint2}, we obtain that
$$
(A^j b)^{\top} \frac{\partial^2h(x)}{\partial x^2} A x+\frac{\partial h(x)}{\partial x} A^{j+1} b-\alpha \frac{\partial h(x)}{\partial x} A^j b =\frac{\partial }{\partial x}\left(\frac{\partial h(x)}{\partial x}Ax-\alpha h(x)\right)A^jb=0.
$$

% in the direction of $A^j b$. The product rule gives
% $$(A^j b)^{\top} \nabla^2 h(x) A x+\nabla h(x)^{\top} A^{j+1} b-\alpha \nabla h(x)^{\top} A^j b=0.$$
The first term on the left-hand side of the equation equals zero by the symmetry of the Hessian matrix and~\eqref{eq:hessian_Ax}. The third term equals zero by~\eqref{eq:induction_assump}. Thus, we have that
$$\frac{\partial h(x)}{\partial x} A^{j+1} b=0.$$
Since
$\operatorname{col}_\R \mathcal C(A,b)$ is spanned
by $b,Ab,\ldots,A^{n-1}b$,~\eqref{eq:cond_controllability}
follows by linearity. This completes our proof.
\end{proof}

%With $A$ being diagonalizable over $\C$, we are able to utilize the projected coordinate 
%$y=Qx$ to retain exactly the state information on which the nonlinear component of an admissible polynomial observable can depend.

For the moment, we assume that $q > 0$. Let $Q$ and $\lambda$ be given as in Lemma~\ref{lem:ker_Q}. For convenience, let $\Lambda:= \diag(\lambda)$. 
Also, recall that $y = Qx$, where $Q\in \C^{q\times n}$ with $q\leq n$. In the following lemma, we translate the condition~\eqref{eq:cond} which is expressed in~$x$ to a new condition expressed in~$y$. 
%The following lemma characterizes the condition for the polynomial on $y$ to satisfy for~\eqref{eq:cond} to hold.

\begin{Lemma}
\label{lem:iff_z}
     A polynomial $h:\R^n\rightarrow\R$ satisfies~\eqref{eq:cond} if and only if there exists a complex polynomial $\hat h$ such that 
    \begin{equation}
    \label{eq:hQx}
        h(x)=\operatorname{Re} \left[\hat h(Qx)\right], \ x\in \R^n,
    \end{equation}
    and 
    \begin{equation}
    \label{eq:invariance_z}
    \sum_{i=1}^q \lambda_i y_i \frac{\partial \hat{h}}{\partial y_i}(y)=\alpha \hat{h}(y), \quad y\in \C^q. 
    \end{equation}
\end{Lemma}
\begin{proof}
We first assume that~\eqref{eq:cond} holds and exhibit a desired~$\hat h$. We extend the domain of~$h$ to $\C^n$; with slight abuse of notation, we still use~$h$ to denote this complex polynomial (with real coefficients).  %it to $\C^n$ using the same polynomial expression
%and denote as $h $. 
% As polynomial
% identities on $\R^n$ remain valid on
% $\C^n$, 
By Lemmas~\ref{lem:ker_Q} and~\ref{lem:iff}, we have that $\frac{\partial h (z)}{\partial z} v=0$ for any $z\in\C^n$ and any $v\in\ker Q$ and hence, 
    \begin{equation}
    \label{eq:const_motion}
        \frac{d}{d t} h (z+t v)=\frac{\partial h (z+t v)}{\partial z} v=0, \quad t\in \R.
    \end{equation}
    which shows that $h $ is constant on any affine subspace of $\C^n$ parallel to $\ker Q$.
    Since $Q$ has full row rank,
its Moore--Penrose inverse is given by
$P:=Q^*(QQ^*)^{-1}\in\C^{n\times q}$,
where $Q^*$ denotes the conjugate
transpose of $Q$. 
    Note that $QP=I_q$, 
    we then define the polynomial 
    \begin{equation}
    \label{eq:def_h_C_hat}
        \hat h(y):=h (Py).
    \end{equation}
    For any $z\in \C^n$, we have that
    $Q(z-PQz)=Qz-QPQz=0$, so $z-PQz \in \ker Q$. 
    It then follows from~\eqref{eq:const_motion} that $h (z)=h (PQz)=\hat h (Qz)$, 
    % In particular, for every $x\in \R^n$, 
    % $$ h(x)=\hat h(Qx),$$
    so $\hat h$ satisfies~\eqref{eq:hQx}. 
    Moreover, we have that
    \begin{equation}
    \label{eq:nabla_hc}
    % \begin{aligned}
        \frac{\partial h(z)}{\partial z} A z=\left.\frac{\partial \hat h(y)}{\partial y} \right |_{y = Q z} Q A z
        =\left.\frac{\partial\hat h(y)}{\partial y} \right |_{y = Q z}  \Lambda Q z, %=\sum_{i=1}^q \lambda_iy_i\frac{\partial \hat h}{\partial y_i}        
    % \end{aligned}
    \end{equation}
    where the first equality follows from the chain rule and $h(z)=\hat h(Qz)$, and the second equality follows from the fact that $QA=\Lambda Q$ as shown in Lemma~\ref{lem:ker_Q}.  
    Combining the above arguments, we have that for any $y\in \C^q$,  
$$
\begin{aligned}
\alpha\hat h(y) &=\alpha h (Py)
=\left.\frac{\partial h(z)}{\partial z}\right |_{z=Py}  APy=\left.\frac{\partial \hat h(z)}{\partial z}\right |_{z=QPy} \Lambda QPy=\frac{\partial \hat h(y)}{\partial y} \Lambda y=\sum_{i=1}^q \lambda_i y_i \frac{\partial \hat{h}}{\partial y_i}(y),
\end{aligned}
$$
where the first equality follows directly from~\eqref{eq:def_h_C_hat}, the second equality follows by the second identity 
in~\eqref{eq:cond} which remains valid over $\C^n$, the third equation follows from~\eqref{eq:nabla_hc} with $z=Py$, the fourth equality follows from $QP=I_q$ and the last equality follows directly from the definition of $\Lambda$. 
Thus, $\hat h$ satisfies~\eqref{eq:invariance_z}.

    Conversely, we assume that $\hat h$ satisfies~\eqref{eq:invariance_z}, %and $\hat h(Qx)$ is real-valued for any $x\in \R^n$, 
    and show that $h(x):= \operatorname{Re}\left[\hat h(Qx)\right]$ satisfies~\eqref{eq:cond}. 
    %For any $x\in \R^n$, 
    Let $y:=Qx$. Then, by the chain rule, 
    $$\frac{\partial h(x)}{\partial x} b=
\operatorname{Re}\left[\left.\frac{\partial\hat h(y)}{\partial y}\right|_{y=Qx} Qb\right]=0,$$
    where the first equality follows from~\eqref{eq:hQx} and the second equality 
    follows from $Qb=0$.
    Moreover, we have that
    \begin{equation*}
        \begin{aligned}
\frac{\partial h(x)}{\partial x} Ax
=
\operatorname{Re}\left[\left.\frac{\partial\hat h(y)}{\partial y}\right|_{y=Qx} QAx\right]
=
\operatorname{Re}\left[\left.\frac{\partial\hat h(y)}{\partial y}\right|_{y=Qx} \Lambda Qx\right]
=\operatorname{Re}\left[\alpha\hat h(Qx)\right]=\alpha h(x),
\end{aligned}
    \end{equation*}
    where the first equality follows from~\eqref{eq:hQx} and the chain rule, the second equality follows from $QA=\Lambda Q$, the third equality follows from~\eqref{eq:invariance_z} with $y=Qx$, and the last equality follows from $\alpha\in\R$ and~\eqref{eq:hQx}. This completes our proof.
\end{proof}

With the above results, we now prove Proposition~\ref{prop:y_alpha}.
% \begin{proposition}
% \label{prop:y_alpha}
%     Let $\hat{h}^{(m)}:\C^q\rightarrow \C$ be homogeneous polynomial of degree $m\geq 1$. Then~\eqref{eq:invariance_z} holds if and only if 
%     \begin{equation}
%     \label{eq:hhat_z}
%         \hat{h}^{(m)}(y)=\sum_{\substack{\beta \in \mathbb{N}_0^q \\|\beta|=m \\ \lambda^{\top} \beta=\alpha}} c_\beta y^\beta
%     \end{equation}
%     where the coefficients $c_\beta\in\C$ are arbitrary.
% \end{proposition}
\begin{proof}[Proof of Proposition~\ref{prop:y_alpha}]
We establish the two items subsequently. 

\xc{Proof of item~1 (case $q = 0$).} 
It is clear that if $h=0$, then~\eqref{eq:cond} holds. 
We now assume that~\eqref{eq:cond} holds and show that $h=0$. 
Since $q = 0$, $ \operatorname{col}_\R \mathcal{C}(A,b)= \R^n$. 
Then, by Lemma~\ref{lem:iff}, $h$ satisfies 
$\frac{\partial h(x)}{\partial x} v=0$, for any $x,v\in \R^n$, so $\frac{\partial h(x)}{\partial x}=0$ for all $x\in \R^n$. Since $h$ has no nonzero constant term, we have that $h=0$.

\xc{Proof of item~2 (case $q>0$).}    
%We now consider the case where $q>0$.
We first assume that~\eqref{eq:cond} holds
and show that $h$ can be written in the form
of~\eqref{eq:h_spectral}. 
By Lemma~\ref{lem:iff_z}, there exists
a complex polynomial $\hat h$ with no nonzero
constant or nonzero linear terms such that~\eqref{eq:hQx} and~\eqref{eq:invariance_z}
hold. We write
\begin{equation}
\label{eq:def_hat_h}
    \hat h(y)=\sum_{|\beta|\geq2}c_\beta y^\beta,
\end{equation}
where the summation is finite and $c_\beta\in \C$. 
For each monomial $y^\beta$, we have that
$y_i\frac{\partial y^\beta}{\partial y_i}
=\beta_i y^\beta$. 
It then follows that
\begin{equation}
\label{eq:lambda_beta}
\sum_{i=1}^q\lambda_i y_i
\frac{\partial y^\beta}{\partial y_i}
=(\lambda^\top\beta)y^\beta.
\end{equation}
Now, substituting~\eqref{eq:def_hat_h}
into~\eqref{eq:invariance_z} and using~\eqref{eq:lambda_beta}, we obtain that
$$
\sum_{i=1}^q \lambda_i y_i \frac{\partial \hat{h}}{\partial y_i}(y)-\alpha \hat{h}(y)= \sum_{|\beta|\geq2}c_\beta
\lambda^\top\beta  y^\beta-\sum_{|\beta|\geq2}c_\beta\alpha  y^\beta=\sum_{|\beta|\geq2}c_\beta (\lambda^\top \beta-\alpha) y^\beta=0,
$$
which holds for all $y\in \C^q$. 
By the linear independence of distinct
monomials over $\C$, we have that
$(\lambda^\top\beta-\alpha)c_\beta=0$, 
for any $\beta$ with $|\beta|\geq2$.
Thus, $c_\beta$ is nonzero only when
$\lambda^\top\beta=\alpha$, i.e., $\beta\in \mathcal S_\alpha^2$. 
% Therefore, we have that
% $\beta \in \mathcal S_\alpha^2$. 
Combining this fact with~\eqref{eq:hQx}, we have that $h$ takes the form~\eqref{eq:h_spectral}.

Conversely, we show that if $h$ is of the form~\eqref{eq:h_spectral}, then it satisfies~\eqref{eq:cond}. 
With the subset $\mathcal{T}\subseteq\mathcal{S}^2_\alpha$ and the coefficients $c_\beta$ given in~\eqref{eq:h_spectral}, we define 
$$
\hat h(y):=\sum_{\beta\in\mathcal T}
c_\beta y^\beta.
$$
Note that any $\beta\in \mathcal{S}^2_\alpha$ satisfies
$|\beta|\geq2$, so $\hat h$ does not contain nonzero constant nor nonzero linear term. 
Moreover, any such $\beta$ satisfies
$\lambda^\top\beta=\alpha$ for any
$\beta\in \mathcal{T}$. Thus, for any $\beta\in \mathcal{T}$ and for any $y\in \C^q$, 
$$
\sum_{i=1}^q\lambda_i y_i
\frac{\partial y^\beta}{\partial y_i}
=\left(\sum_{i=1}^q\lambda_i\beta_i\right)y^\beta
=(\lambda^\top\beta)y^\beta
=\alpha y^\beta.
$$
It follows that $\hat h$ satisfies~\eqref{eq:invariance_z}. 
Also, since~$h$ is of the form~\eqref{eq:h_spectral}, we have that~\eqref{eq:hQx} holds. 
By Lemma~\ref{lem:iff_z}, we conclude that $h$ satisfies~\eqref{eq:cond}. 
This completes our proof.
\end{proof}

\subsection{Step 2: Connection with finite generators}

We now show that the representation~\eqref{eq:h_spectral} can be expressed
in the generators of $\mathcal S_0$ and
the minimal elements of
$\mathcal S_\alpha^2$. Specifically, we have the following proposition.
% \begin{Lemma}
% \label{lem:minimal_elements}
%     Let $\mathcal S_\alpha^{\geq 2}$ be defined as in~\eqref{eq:semigroup}, and its minimal elements be
%     $\kappa^{(1)},\ldots,\kappa^{(N)}$.
%     Then
%     $$\mathcal{S}_\alpha^{\geq 2}=\bigcup_{j=1}^N\left(\kappa^{(j)}+\mathcal{S}_0\right).$$
% \end{Lemma}
% \begin{proof}
%     by Dickson's lemma, $\mathcal{S}_\alpha^{\geq 2}\subseteq \N_0^q$ has finitely many componentwise-minimal elements, denoted by $\kappa^{(1)},\ldots,\kappa^{(N)}$.

%     Consider any $\beta\in \mathcal{S}_{\alpha}^{\geq 2}$. There exists a minimal element $\kappa^{(j)}\leq\beta$. We define $\eta:=\beta-\kappa^{(j)}$. By the definition of $\mathcal{S}_{\alpha}^{\geq 2}$, we have that 
%     $$\lambda^\top \eta=\lambda^\top \beta-\lambda^\top \kappa^{(j)}=0.$$ Hence $\eta\in \mathcal{S}_0$ and $\beta\in \kappa^{(j)}+\mathcal{S}_0$.

%     Conversely, if $\eta\in \mathcal{S}_0$, then $\kappa^{(j)}+\eta$ is nonnegative and it satisfies that 
%     $$\lambda^\top (\kappa^{(j)}+\eta)=\alpha,$$
%     where $\left|\kappa^{(j)}+\eta\right| \geq\left|\kappa^{(j)}\right| \geq 2$. Therefore, $\kappa^{(j)}+\eta \in \mathcal{S}_{\alpha}^{\geq 2}$.
%     This completes our proof.
% \end{proof}

% {\color{blue}
% bring back the proposition; and we take the two steps to represent $h$;

% reorganize section 3 as proof of theorem 1
% }

\begin{proposition}
\label{prop:finite}
Let $\eta^{(1)},\ldots,\eta^{(K)}$ be
generators of $\mathcal S_0$, and
$\kappa^{(1)},\ldots,\kappa^{(N)}$ be
the minimal elements of
$\mathcal S_\alpha^2$.
There exist a finite subset $\mathcal{T}\subseteq \mathcal{S}_\alpha^2$ and coefficients $\{c_\beta  \mid \beta \in \mathcal{T}\}$ such that
\begin{equation}
\label{eq:hhat_spectral}
\hat h(y)=
\sum_{\beta\in\mathcal T}
c_\beta y^\beta, \quad y\in\C^q, 
\end{equation}
if and only if
\begin{equation}
\label{eq:hhat_repr}
\hat h(y)=
\sum_{j=1}^N y^{\kappa^{(j)}}
F_j\left(
y^{\eta^{(1)}},\ldots,
y^{\eta^{(K)}}
\right),
\end{equation}
where $F_1,\ldots,F_N$ are arbitrary
complex polynomials.
\end{proposition}

\begin{proof}
    We first assume that~\eqref{eq:hhat_spectral}
holds and prove that $\hat h$ can be written
in the form of~\eqref{eq:hhat_repr}.
For any $\beta\in \mathcal{T}$, there exists a
minimal element
$\kappa^{(j)}\leq\beta$.
Let $\eta:=\beta-\kappa^{(j)}$.
Since both $\beta$ and $\kappa^{(j)}$
belong to $\mathcal S_\alpha^2$, we have that $\eta\in\N_0^q$ and
$$
\lambda^\top\eta
=\lambda^\top\beta
-\lambda^\top\kappa^{(j)}=0,
$$
and hence, $\eta\in\mathcal S_0$.
By Lemma~\ref{lem:finite_generators},
there exist $a_1,\ldots,a_K\in\N_0$
such that
$\eta=\sum_{i=1}^K a_i\eta^{(i)}.$
It follows that
\begin{equation*}
y^\beta=y^{\kappa^{(j)}+\eta}=y^{\kappa^{(j)}+\sum_{i=1}^K a_i\eta^{(i)}}
=y^{\kappa^{(j)}}\prod_{i=1}^K
\left(y^{\eta^{(i)}}\right)^{a_i}.
\end{equation*}
Substituting the above into~\eqref{eq:hhat_spectral} and grouping the monomials
of $\hat h$ associated with each
$\kappa^{(j)}$ gives~\eqref{eq:hhat_repr}.

Conversely, we show that if $\hat h$ is of the
form~\eqref{eq:hhat_repr}, then it satisfies~\eqref{eq:hhat_spectral}. We expand
each $F_j$ into its finitely many
monomials and represent the exponent
of each term as
$\beta=\kappa^{(j)}
+\sum_{i=1}^K a_i\eta^{(i)}$, with $a_i\in\N_0$.
It follows that
$$
\lambda^\top\beta
=\lambda^\top\kappa^{(j)}
+\sum_{i=1}^K a_i
\lambda^\top\eta^{(i)}
=\alpha.
$$
By the definition of $\kappa^{(j)}$, we have that
$|\beta|\geq|\kappa^{(j)}|\geq2$. Thus, $\beta\in\mathcal S_\alpha^2$, which implies that every monomial in~\eqref{eq:hhat_repr} has an exponent
in $\mathcal S_\alpha^2$.
Since each $F_j$ is a complex polynomial,
$\hat h$ contains finitely many
monomials and therefore admits
the representation~\eqref{eq:hhat_spectral}.
This completes our proof.
\end{proof}

Theorem~\ref{thm:main} then follows directly from
Lemma~\ref{lem:matching} and
Propositions~\ref{prop:y_alpha} and~\ref{prop:finite}.
By Lemma~\ref{lem:matching},
$p$ satisfies~\eqref{eq:invariance}
if and only if $(e,\ell)$ satisfies
~\eqref{eq:constraints} and
$h$ satisfies~\eqref{eq:cond}.
For $q=0$, Proposition~\ref{prop:y_alpha}
shows that~\eqref{eq:cond} holds
if and only if $h=0$.
For $q>0$, Proposition~\ref{prop:y_alpha}
shows that $h$ satisfies~\eqref{eq:cond}
if and only if it takes the form~\eqref{eq:h_spectral}.
It then follows from
Proposition~\ref{prop:finite},
by substituting $y=Qx$ and taking
real parts, that $h$ admits
the representation~\eqref{eq:h_form},
and conversely. Thus, $p$ satisfies~\eqref{eq:invariance} if and only if the conditions in Theorem~\ref{thm:main} hold.
%where $\hat h$ admits a
%representation of~\eqref{eq:hhat_spectral}.
This completes our proof. \qed

% With the results established in Section~\ref{sec:compatibility}, we prove Theorem~\ref{thm:main} in this section. It remains to show that there exist a finite subset $\mathcal{T}\subseteq \mathcal{S}_\alpha^2$ and coefficients $\{c_\beta  \mid \beta \in \mathcal{T}\}$ such that
% \begin{equation}
% \label{eq:hhat_spectral}
% \hat h(y)=
% \sum_{\beta\in\mathcal T}
% c_\beta y^\beta,
% \end{equation} if and only if $\hat h(y)$ is in the form of 
% \begin{equation}
% \label{eq:hhat_repr}
% \hat h(y)=
% \sum_{j=1}^N y^{\kappa^{(j)}}
% F_j\left(
% y^{\eta^{(1)}},\ldots,
% y^{\eta^{(K)}}
% \right),
% \end{equation}
% where $\eta^{(1)},\ldots,\eta^{(K)}$ are
% generators of $\mathcal S_0$,
% $\kappa^{(1)},\ldots,\kappa^{(N)}$ are
% the componentwise minimal elements of
% $\mathcal S_\alpha^2$, and $F_1,\ldots,F_N$ are arbitrary
% complex polynomials.

\section{Conclusion and Outlook}
\label{sec:conclusion}
In this paper, we introduced and studied
an inverse problem for super-linearization, which seeks to characterize all dynamical systems that can be strongly super-linearized to a prescribed linear system. Focusing on the case where the dimension reduction is one and the leading $n$-by-$n$ submatrix is
diagonalizable, we provide a complete characterization of all compatible real polynomial
observables.
% strong superlinearization with
% a single polynomial observable $p$ for a
% prescribed lifted matrix $M$, whose submatrix
% $A$ is diagonalizable over $\C$.
% We derived necessary and sufficient conditions
% on the affine and nonlinear components of $p$
% for $(M,p)$ to be a strong superlinearization
% of the corresponding nonlinear system.
% These conditions provide a finite-generator
% representation using the generators of
% $\mathcal S_0$ and the minimal
% elements of $\mathcal S_\alpha^2$. The results provide a complete
% characterization of all real polynomial
% observables $p$ that is compatible with $M$ across all
% admissible degrees.
In our future work, we will first extend the present
characterization to general  matrices $M$ (by relaxing the diagonalizability assumption on the leading principal submatrix $A$) and
then to include multiple
polynomial observables.

%\nocite{*}
\printbibliography

\appendix

{

% \section{Proof of Lemma~\ref{lem:ker_Q}}
% \begin{proof}
    
% \end{proof}

\section*{Appendix}

\section{Proof of Lemma~\ref{lem:matching}}
\label{app:matching_proofs}

\begin{proof}
    We first assume that $p$ satisfies~\eqref{eq:invariance}, and show that $(e,\ell)$ satisfies~\eqref{eq:constraints} and $(h,\alpha)$ satisfies~\eqref{eq:cond}. 
    We prove the result by matching homogeneous terms of the same degree on both sides of~\eqref{eq:invariance}. If $h=0$, then~\eqref{eq:constraints} and~\eqref{eq:cond} directly follow. We therefore assume that $h\neq 0$ and write $h(x)=\sum_{m=2}^k h^{(m)}(x)$ where $k\geq 2$, and $h^{(m)}$ denotes the homogeneous polynomial of degree $m$.
    
    Substituting~\eqref{eq:px}  into~\eqref{eq:invariance}, by computation, we have that
    \begin{align}
    %\label{eq:degree_match}
        &\frac{\partial p(x)}{\partial x}
(Ax+bp(x)) \notag\\
={}&
\left(\ell^\top+
\sum_{m=2}^k
\frac{\partial h^{(m)}(x)}{\partial x}\right)
\left[
Ax+b\left(e+\ell^\top x+
\sum_{j=2}^k h^{(j)}(x)\right)
\right]\notag\\
={}&
\ell^\top(Ax+b\ell^\top x)
+\ell^\top be
+\ell^\top b\sum_{j=2}^k h^{(j)}(x)
+\sum_{m=2}^k
\frac{\partial h^{(m)}(x)}{\partial x}Ax
+e\sum_{m=2}^k
\frac{\partial h^{(m)}(x)}{\partial x}b \notag\\
&+\ell^\top x\sum_{m=2}^k
\frac{\partial h^{(m)}(x)}{\partial x}b+\sum_{m=2}^k
\frac{\partial h^{(m)}(x)}{\partial x}b
\sum_{j=2}^k h^{(j)}(x)\label{eq:degree_match}\\
={}&c^\top x+
d\left(e+\ell^\top x+
\sum_{m=2}^k h^{(m)}(x)\right).\notag
\end{align}

    In the following, we group the terms with the same degree.
    First, notice that the only term of the highest degree $(2k-1)$ in~\eqref{eq:degree_match} is $\frac{\partial h^{(k)}(x)}{\partial x}
bh^{(k)}(x)$. Since $h^{(k)}$ is a nonzero polynomial, it follows that \begin{equation}
    \label{eq:pde}
        \frac{\partial h^{(k)}(x)}{\partial x}b=0.
    \end{equation} 
    For $k\geq 3$, we next consider the terms of degree $2k-2$. To match the degree in equation~\eqref{eq:degree_match}, necessarily we have 
    $$\begin{aligned}
&\frac{\partial h^{(k)}(x)}{\partial x}
b\,h^{(k-1)}(x)+\frac{\partial h^{(k-1)}(x)}{\partial x}
b\,h^{(k)}(x)=0.
\end{aligned}$$
    By~\eqref{eq:pde}, we obtain
    \begin{equation}
        \frac{\partial h^{(k-1)}(x)}{\partial x}b=0.
    \end{equation}
    We now prove by descending induction
    that $\frac{\partial h^{(s)}(x)}{\partial x}b=0$ holds for any $2\leq s< k-1$.
    Assume that, for some $2\leq s< k-1$, we have 
    \begin{equation}
        \label{eq:induction}
         \frac{\partial h^{(k-1)}(x)}{\partial x} b= \frac{\partial h^{(k-2)}(x)}{\partial x} b=\cdots = \frac{\partial h^{(s+1)}(x)}{\partial x} b=0.
    \end{equation}
    We next show that $\frac{\partial h^{(s)}(x)}{\partial x} b=0$. 
    Consider the terms with degree $s+k-1$ in~\eqref{eq:degree_match}, the corresponding terms are found as
    $$\sum_{\substack{2\leq i,j\leq k\\
i+j-1=s+k-1}}
\frac{\partial h^{(i)}(x)}{\partial x}
b\,h^{(j)}(x)=0.$$ By~\eqref{eq:degree_match} and~\eqref{eq:induction}, the remaining component is 
    $\frac{\partial h^{(s)}(x)}{\partial x} bh^{(k)}(x)$. Since the highest degree polynomial $h^{(k)}(x)$ is nonzero, we have that $\frac{\partial h^{(s)}(x)}{\partial x} b=0$. This establishes~\eqref{eq:constraint1}.

    We then consider the homogeneous terms in~\eqref{eq:degree_match} of degree $j$ with $2\leq j\leq k$. After applying~\eqref{eq:constraint1}, the remaining homogeneous terms of degree $j$ satisfy:
    \begin{equation}
    \label{eq:h_identity}
        \frac{\partial h^{(j)}(x)}{\partial x} Ax=\alpha h^{(j)}(x).
    \end{equation}
    Summing~\eqref{eq:h_identity} over $j=2,\ldots,k$ proves~\eqref{eq:constraint2}.

    It remains to show that~\eqref{eq:constraints} holds. Substitute~\eqref{eq:constraint1} into~\eqref{eq:invariance}, for any $x\in \R^n$, the remaining linear term satisfies 
    $$(\ell^\top A+\ell^\top b\ell^\top -c^\top -d\ell^\top)x =0.$$
    Rearrange the coefficients to yield that $$c^\top-\ell^\top A+(d-\ell^\top b)\ell^\top=c^\top-\ell^\top A+\alpha\ell^\top=0.$$
    Then~\eqref{eq:constraint4} follows directly by taking the transpose. Finally, we compare the constant
terms in~\eqref{eq:degree_match}, and obtain
    $\alpha e=0$. This establishes~\eqref{eq:constraint3}.

    Conversely, we assume that~\eqref{eq:constraints} and~\eqref{eq:cond} hold and show that $p$ satisfies~\eqref{eq:invariance}. Substituting~\eqref{eq:px}
into the left-hand side of~\eqref{eq:invariance}, we have that
$$
\begin{aligned}
\frac{\partial p(x)}{\partial x}(Ax+bp(x))
&=\ell^\top Ax+\frac{\partial h(x)}{\partial x} Ax
+\ell^\top b\,p(x)
+\frac{\partial h(x)}{\partial x} b\,p(x)\\
&=\ell^\top Ax+\alpha h(x)
+\ell^\top b\,p(x)\\
&=(\ell^\top A-\alpha\ell^\top)x
+(\alpha+\ell^\top b)p(x)-\alpha e\\
&=c^\top x+dp(x),
\end{aligned}
$$
where the second equality follows from~\eqref{eq:cond}, and the last equality
follows from~\eqref{eq:constraints}
and $\alpha=d-\ell^\top b$.
Thus,~\eqref{eq:invariance} holds.
This completes our proof.
\end{proof}

\end{document}